\documentclass[12pt]{article}
\newtheorem{theorem}{Theorem}

\newtheorem{lemma}[theorem]{Lemma}
\newtheorem{corollary}[theorem]{Corollary}
\newcommand\QED{\rule{1.5mm}{1.5mm}}
\newenvironment{proof}{\par\noindent{\bf Proof:}\rm\enspace}{\QED\par}
\begin{document}

\title{Iterating Random Functions on a Finite Set}
\author{William M.Y. Goh\\
\small MCS Department, Drexel University\\
\small Philadelphia, Pa. 19104 \\
\small wgoh@mcs.drexel.edu \and
Pawel Hitczenko\\
\small MCS Department, Drexel University\\
\small Philadelphia, Pa. 19104 \\
\small phitczen@mcs.drexel.edu \and Eric
Schmutz\\\small  MCS Department, Drexel University\\
\small Philadelphia, Pa., 19104\\
 \small Eric.Jonathan.Schmutz@drexel.edu\\
 }

\maketitle \abstract{
 Choose random functions $f_1,
f_2,f_3,\dots $ independently and uniformly from among the $n^n$
functions from $[n]$ into $[n]$. For $t>1$, let $g_t=f_t\circ
f_{t-1}\circ \cdots \circ f_1$ be the composition of the first $t$
functions, and let  $T$ be the smallest $t$ for which $g_t$ is
constant(i.e. $g_t(i)=g_t(j)$ for all $i,j$).  We prove that, for
any positive real number $x$,
$$\lim\limits_{n\rightarrow\infty}\Pr({T\over n}\leq x)=
\int\limits_{0}^{x}f(y)dy,$$ where
$$f(y)=\sum\limits_{k\geq 2}(-1)^{k}e^{-y{k\choose
2}}(2k-1){k\choose 2}.$$ } We make our proof available here, but
acknowledge that the the result is already \lq\lq well
known.\rq\rq

\par

 \vfill\eject
\section{Introduction}

 Let
$f_1,f_2,f_3,\dots $ be a sequence of  functions chosen
independently and uniformly randomly from the $n^n$ functions on
$[n]$. Let $g_1 = f_{1}$, and for $t>1$ let $g_{t}=f_{t}\circ
g_{t-1}$ be the composition of the first $t$ random functions.
Define $T(\langle f_{i}\rangle_{i=1}^{\infty})$ to be the smallest
$t$ for which $g_t$ is a constant function. (i.e. $g_t(i)= g_t(j)$
for all $i\not=j.$) This manuscript contains a simple derivation
of the asymptotic distribution of $T$.
  We had originally intended to publish it
in a journal, but  we recently learned that the the asymptotic
distribution of $T$
 is \lq\lq well
known\rq\rq.  It was apparently  known to Kingman twenty years ago
\cite{King},\cite{King2},\cite{King3}, and stronger results are
fully proved in Donnelly\cite{Don}.  There is a lot of  related
work by M\H{o}hle and others,e.g. \cite{mohle1},\cite{mohle2}, and
\cite{Fill}.

\par
For $m>1$, let $\tau_{m}=\bigl|\bigl\lbrace t:
|Range(g_t)|=m\bigr\rbrace \bigr| $ be the the number of iterates
for which the range has exactly $m$ elements. Thus
  $T= \sum_{m=2}^{n}\tau_{m}$.
  The random variables $\lbrace \tau_{m}\rbrace_{m=2}^{n}$ are  not
independent.
They are however {\sl conditionally } independent once we specify
the set of visited states.  Fortunately this set  is well behaved,
has some convenient properties that enable us to do computations.

Let $\xi=\lfloor \log\log n\rfloor$,
 and decompose $T$
as $T=T_1+T_2$, where $T_1=\sum\limits_{m=2}^{\xi} \tau_m$ and
$T_2=\sum\limits_{m=\xi+1}^n\tau_m$.
 Let
 ${\cal A}=\bigcap\limits_{m=1}^{\xi}[\tau_m > 0]$. The
following facts from \cite{Avi} will be needed (See also Theorem 5
of  \cite{King3}):
\begin{theorem}
\label{facts} $\Pr({\cal A})=1-o(1),$ and $E(T_{2})=o(n).$
\end{theorem}

\section{Characteristic Function}
  Let
$\lambda_{k}=\prod\limits_{j=1}^{k-1}(1-{j\over n}).$ Then we have
\label{Char}
\begin{theorem}
$E(e^{itT_{1}}|{\cal A}
)=e^{it(\xi-1)}\prod\limits_{k=2}^{\xi}{(1-\lambda_{k})\over
1-\lambda_{k}e^{it}}.$
\end{theorem}
\begin{proof}

Suppose $g_{t-1}$ has an $m$ element range $R=\lbrace
r_1,r_2,\dots ,r_m\rbrace $. What is the chance that the next
function $g_t$ still has an $m$ element range?
  On $R$ we have $n$
choices for $f_t(r_1)$,\ then $n-1$ choices for $f_t(r_2 )$ etc.
For $x\notin R$, $f_t(x)$ can be chosen arbitrarily. Hence the
number of functions $f_t$ for which $g_t=f_t\circ g_{t-1}$ has an
an $m$ element range is $n^{n-m}\prod\limits_{j=0}^{m-1}(n-j).$
Hence
\begin{equation} \Pr(\tau_m =k|\tau_m >0)=
\lambda_{m}^{k-1}(1-\lambda_m),
\end{equation} and consequently
$$
E(e^{it\tau_m}|\tau_m>0)=
\sum\limits_{k=1}^{\infty}\lambda_{m}^{k-1}(1-\lambda_{m})e^{ikt}=
{(1-\lambda_{m})e^{it}\over 1-\lambda_{m}e^{it}}
$$
\end{proof}
Now let $\phi_{n}(t)=E(e^{itT_{1}/n}|{\cal A})$ be the
characteristic function of the normalized random variable
$T_{1}/n$ on ${\cal A}.$ Then the following corollary follows
immediately from  Theorem \ref{Char}.
\begin{corollary}
$\phi_{n}(t)=\prod\limits_{m=2}^{\xi} {
(1-\lambda_{m})e^{it/n}\over 1-\lambda_{m}e^{it/n}}=e^{it(\xi
-1)/n} \prod\limits_{m=2}^{\xi} { (1-\lambda_{m})\over
1-\lambda_{m}e^{it/n}}
 $
\end{corollary}

\begin{lemma}
 $\phi_{n}(t)=
\prod\limits_{m=2}^{\infty}{{m\choose 2}\over {m\choose
2}-it}+o(1).$
\end{lemma}
\begin{proof}
Note that, for $m\leq \xi$,
\begin{equation} 1-\lambda_{m}={1\over n}{m\choose
2}+O({\xi^{4}\over n^{2}})
\end{equation}
and
\begin{equation} 1-\lambda_{m}e^{it/n}={1\over n}\Bigl({m\choose
2}-it\Bigr)+O({\xi^{4}\over n^{2}})
\end{equation}
Therefore
\begin{equation}
 { (1-\lambda_{m})e^{it/n}\over
1-\lambda_{m}e^{it/n}}={{{m\choose 2}+O({\xi^{4}\over n})}\over
{m\choose 2}-it+O({\xi^{4}\over n})} ={{{m\choose 2}}\over
{m\choose 2}-it}\bigl(1+O({\xi^{4}\over n})\bigr).
\end{equation}
Therefore $$\phi_{n}(t)=\bigl(1+O({\xi^{4}\over n})\bigr)^{\xi}
\prod\limits_{m=2}^{\xi}{{{m\choose 2}}\over {m\choose 2}-it}$$
$$=(1+o(1))\prod\limits_{m=2}^{\xi}{{{m\choose 2}}\over {m\choose 2}-it}.$$
Finally, note that the infinite product
$\prod\limits_{m=2}^{\infty}{{{m\choose 2}}\over {m\choose 2}-it}
=\prod\limits_{m=2}^{\infty}{1\over 1-{it\over {m\choose 2}}} $
converges since $\sum\limits {m\choose 2}^{-1}$ is convergent.

\end{proof}

\section{Simplification}
\label{Nice} To facilitate inversion, we reexpress the
characteristic function $\phi$ in a more convenient form. Working
with the reciprocal, we have
\begin{equation}
\label{recip} {1\over \phi_{n}(t)+o(1)}=\prod\limits_{k\geq
1}\bigl(1-{2it\over k(k+1)}\bigr) =\prod\limits_{k\geq
1}{(k-\alpha)(k-\beta)\over k(k+1)}, \end{equation}
 where
$\alpha={{-1-\sqrt{1+8it}}\over 2}, \beta={{-1+\sqrt{1+8it}}\over
2}.$ It is well known \cite{AS} that $${1\over
\Gamma(z)}=ze^{\gamma z}\prod\limits_{n\geq 1}\Bigl((1+{z\over
n})e^{-z/n}\Bigr).$$ Since $\alpha+\beta=-1$, the right side of
equation (\ref{recip}) becomes
$$ \prod\limits_{k\geq
1}{(1-\alpha/k)e^{\alpha/k}(1-\beta/k)e^{\beta/k} \over (1+{1\over
k})e ^{-1/k}}= {1\over \alpha\beta
\Gamma(-\alpha)\Gamma(-\beta)}={\cos({\pi\over
2}\sqrt{1+8it})\over -2\pi it}.$$ Hence
\begin{equation}
\label{nicechar} \phi_{n}(t)= {-2\pi i t\over \cos({\pi\over
2}\sqrt{1+8it})}+o(1).
\end{equation}
\section{Fourier Inversion}
\label{Finis} Inverting, we get the conditional density function
$f_n$ for $T_1/n$:
\begin{equation}\label{invert}
f_{n}(x)={1\over 2\pi}\int\limits_{-\infty}^{\infty}
e^{-itx}({-2\pi i t\over \cos({\pi\over 2}\sqrt{1+8it})}+o(1))dt.
\end{equation}
 Note that ${-2\pi i t\over \cos({\pi\over 2}\sqrt{1+8it})}$
has simple poles at $t=-i{k\choose 2},$ for $k=2,3,4,\dots$ Since
the residue at  $t=-i{k\choose 2}$
is 
$i(-1)^{k}(2k-1){k \choose 2}e^{-{k\choose 2}x}$, contour
integration yields, for $x>0,$
 $f_n(x)=f(x)+o(1)$ where
$$f(x)= \sum\limits_{k\geq 2}(-1)^{k}e^{-{k\choose 2}x}{k\choose
 2}(2k-1).$$
\section{Main Result}
For $x>0$, let $F(x)=\int\limits_{0}^{x}f(t)dt.$
 Our main result is
\begin{theorem}
\label{main} For any $x>0$,
$\lim\limits_{n\rightarrow\infty}\Pr(T/n \leq x)=F(x).$
\end{theorem}
\begin{proof}
 For any $x$,
  $$\Pr(T/n \leq x)\leq \Pr(T_{1}/n\leq x)$$
  $$= \Pr(T_{1}/n\leq x |{\cal A})\Pr({\cal A})+ \Pr(T_{1}/n\leq x |{\cal A}^{c})\Pr({\cal
  A}^{c})$$
  $$=\Pr(T_{1}/n\leq x |{\cal A})(1+o(1))+o(1)$$
  $$=F(x)(1+o(1))+o(1).$$
  In the other direction, let $\epsilon$ be a fixed but
  arbitrarily small positive number. Then
  $$\Pr(T/n \leq x)\geq \Pr(T_1/n \leq x-\epsilon\ {\rm and}\ T_{2}/n
  \leq \epsilon)$$
  $$\geq \Pr(T_{1}/n \leq x-\epsilon)-\Pr(T_{2}/n >\epsilon)$$
 $$ \geq \Pr(T_{1}/n \leq x-\epsilon)- {E(T_{2}/n)\over \epsilon}$$
$$=F(x-\epsilon)+o(1).$$
The theorem follows from this and the fact that $F$ is continuous.
  \end{proof}
\section{Discussion}
Although the ultimate behaviour of our chain is like Kingman's
coalescent \cite{King}, there are differences. In that process
every state is visited, whereas in our process few of the high
numbered states are visited.
 Let
$N=\sum\limits_{m=2}^{n}I_{[\tau_{m}>0]},$ the number of states
visited. In an earlier version of this manuscript, we conjectured
that $E(N)\sim \sqrt{2\pi n}.$ Robin Pemantle recently proved our
conjecture  and the corresponding  central limit theorem. He may
also be able to prove  stronger and more general versions of this
result, e.g. a functional limit theorem. \cite{RP}.

\vfill\eject \noindent{\bf Acknowledgement} There is a large
literature in applied probability that can be traced back to
Kingman's work, and it was not immediately obvious to us what is
relevant.  We are grateful to Simon Tavare for pointing our way to
the the work of Donnelly and  M\H{o}hle.

  \end{document}